\documentclass[11pt,reqno,oneside]{amsart}
\def\R{\mathbb R}
\def\Box{\setlength{\unitlength}{0.01cm}
     \begin{picture}(15,15)(-5,-5)
       \framebox(15,15){}
     \end{picture} }

\newtheorem{theorem}{Theorem}[section]
\newtheorem{lemma}{Lemma}[section]

\newtheorem{remark}{Remark}[section]

\def\proof{\noindent{\bf Proof.}\,\,}
\def\R{{\rm I}\!{\rm R}}
\def\O{\Omega}
\def\uu{{\bf u}}
\def\vv{{\bf v}}
\def\ww{{\bf w}}
\def\nn{{\bf n}}
\def\ve{\varepsilon}
\def\d{\mbox{div\,}}
\def\eps{\varepsilon}

\def\tQ{\widetilde Q}

\begin{document}

\title[]
{An elementary proof of the continuity from $L_0^2(\O)$ to $H^1_0(\O)^n$ of Bogovskii's right inverse
of the divergence}

\author[R. G. Dur\'an]{Ricardo G. Dur\'an}
\address{Departamento de Matem\'atica\\Facultad de Ciencias Exactas y
Naturales\\
Universidad de Buenos Aires and IMAS, CONICET\\Argentina.} \email{rduran@dm.uba.ar}
\thanks{This research was supported by Universidad de Buenos Aires
under grant X070 , by ANPCyT under grant PICT 01307 and by
CONICET under grant PIP 11220090100625.}


\begin{abstract}
The existence of right inverses of the divergence as an operator
form $H^1_0(\O)^n$ to $L_0^2(\O)$ is a problem that has been widely
studied because of its importance in the analysis of the classic equations
of fluid dynamics. When $\O$ is a bounded domain which is star-shaped
with respect to a ball $B$, a right inverse given by an integral operator
was introduced by Bogovskii, who also proved the continuity using
the Calder\'on-Zygmund theory of singular integrals.

In this paper we give an alternative elementary proof using the
Fourier transform. As a consequence, we obtain estimates of the constant
in the continuity in terms of the ratio between the diameters of $\O$ and $B$.
Moreover, using the relation between the existence of
right inverses of the divergence with the Korn and improved Poincar\'e
inequalities, we obtain estimates for the constants in these two inequalities.
\end{abstract}
\maketitle

\section{Introduction}
\label{section0}
\setcounter{equation}{0}

Let $\O\subset\R^n$ be a bounded domain. Given a smooth vector field $\uu$
defined in $\O$ we will denote with $D\uu$ its differential matrix, namely,
$$
D\uu=\left(\frac{\partial u_i}{\partial x_j}\right)
$$
and for a tensor field $(a_{ij})$ we define its norm by
$$
\|a\|^2_{L^2(\O)}=\sum_{i,j=1}^n \|a_{ij}\|^2_{L^2(\O)}.
$$

The existence of solutions $\uu\in H^1_0(\O)^n$ of

\begin{equation}
\label{divu=f}
\d\uu=f
\end{equation}
satisfying

\begin{equation}
\label{estimate}
\left\|D\uu\right\|_{L^2(\O)}
\le C_{div,\O} \|f\|_{L^2(\O)},
\end{equation}
where $f\in L^2(\O)$ has vanishing
mean value, and the constant $C_{div,\O}$ depends only on $\O$, is a problem that
has been widely analyzed because of its several applications and
connections with other important results.

Assume that $\O\subset\R^n$ is a domain with diameter $R$ which is
star-shaped with respect to a ball $B\subset\O$, which we assume centered
at the origin and of radius $\rho$. For a function
$\omega\in C_0^\infty(B)$ such that $\int_\O\omega=1$, a solution
of (\ref{divu=f}) is given by

\begin{equation}
\label{solucion}
\uu(x)= \int_\O G(x,y)f(y)\, dy
\end{equation}
where $G=(G_1,\cdots,G_n)$ is defined by

$$
G(x,y)=\int_0^1 \frac{(x-y)}{t}
\omega\left(y+\frac{x-y}{t}\right)\,\frac{dt}{t^{n}} \ \ .
$$
Moreover, $\uu\in H^1_0(\O)^n$ and (\ref{estimate}) is satisfied.

This formula was introduced in \cite{B} by Bogovskii who proved
the estimate (\ref{estimate}), as well as its generalization for $L^p$, $1<p<\infty$,
using the general Calder\'on-Zygmund theory
of singular integrals developed in \cite{CZ}.

More recently, several papers have considered extensions and
applications of this formula. In \cite{DM}, a weighted version of
(\ref{estimate}), which is of interest in finite element analysis, was proved.
In \cite{ADM}, an extension of Bogovskii's formula was introduced
for the rather general class of John domains and the estimate (\ref{estimate})
was proved using again the Calder\'on-Zygmund theory. Also, extensions of
(\ref{estimate}) for fractional order positive and negative Sobolev norms
have been obtained in \cite{CM,GHH}.

The goal of this paper is twofold:

First, we want to give
a simple proof of the estimate (\ref{estimate}) for the solution given
by (\ref{solucion}) using elementary properties of the Fourier transform.
In this way we avoid the use of the complicated general theory of singular
integral operators. We believe that this can be of interest for
teaching purposes.

Second, we are interested in obtaining some information
on the constant in terms of the ratio $R/\rho$. As a byproduct,
this result can be used to give estimates for the constants
in some Korn and improved Poincar\'e inequalities.

The paper is organized in such a way that the reader interested
only in the first part needs to read only up to the end of Section \ref{section1}
which deals with the continuity of the singular integral operator.
In Section \ref{section2} we modify the proof of the continuity
in order to obtain a sharper estimate of the constant in (\ref{estimate}). Finally, in
Sections \ref{section3} and \ref{section4} we obtain estimates for the
constants in the Korn and improved Poincar\'e inequalities respectively.

\section{Boundedness of the singular integral operator}
\label{section1}
\setcounter{equation}{0}

In order to work with functions defined in $\R^n$ we extend $f$ by
zero outside of $\O$ in (\ref{solucion}).

Let us recall the basic properties of the Fourier
transform that we will need (see for example \cite{S}).
The Fourier transform is defined
for $f\in L^1(\R^n)$ by

$$
\widehat f(\xi)=\int e^{-2\pi i x\cdot\xi} f(x) dx.
$$
Here and in the rest of the paper, when we do not indicate the domain
of integration it is understood that it is $\R^n$. The Fourier transform
can be extended to $f$ in the class of tempered
distributions ${\mathcal S}'$, in particular, it
is defined in $L^2(\R^n)$ and it is an isometry, i. e.,

$$
\|\widehat f\|_{L^2(\R^n)}=\|f\|_{L^2(\R^n)}.
$$
We will use the well known relation

$$
\frac{\widehat{\partial f}}{\partial x_j}(\xi)=2\pi i\xi_j \widehat f(\xi).
$$
The $k$-component of the vector field $\uu$ defined in (\ref{solucion}) is given by

$$
u_k=u_{k,1}-u_{k,2},
$$
where

$$
u_{k,1}(x)=\int_0^1 \int \left(y_k+\frac{(x_k-y_k)}{t}\right)
\omega\left(y+\frac{x-y}{t}\right)\, f(y)\, dy\,\frac{dt}{t^{n}}.
$$
and
$$
u_{k,2}(x)=\int_0^1 \int y_k
\omega\left(y+\frac{x-y}{t}\right)\, f(y)\, dy\,\frac{dt}{t^{n}}.
$$
These double integrals exist, if for example we assume that
$f\in L^1(\R^n)$ and has compact support. Indeed, if $\mbox{supp\,}f\subset B(0,M)$
then, both integrands vanish unless $\left|y+\frac{x-y}{t}\right|<\rho$
and $|y|<M$, and so, assuming that $\rho<M$, we can restrict the domain
of integration to $|x-y|<2Mt$. Therefore, integrating first in the $t$ variable,
it follows that, for $i=1,2$,
$$
|u_{k,i}(x)|\le C \int \frac{|f(y)|}{|x-y|^{n-1}}\, dy,
$$
where the constant $C$ depends only on $\omega$, $n$, and $M$.
Since $f\in L^1(\R^n)$ the last integral is finite for almost every $x$.

In order to take the derivatives of $u_{k,i}$
it is convenient to write

$$
u_{k,1}(x)=\lim_{\eps\to 0}\int_\eps^1 \int \left(y_k+\frac{(x_k-y_k)}{t}\right)
\omega\left(y+\frac{x-y}{t}\right)\, f(y)\, dy\,\frac{dt}{t^{n}},
$$
and
$$
u_{k,2}(x)=\lim_{\eps\to 0}\int_\eps^1\int y_k
\omega\left(y+\frac{x-y}{t}\right)\, f(y)\, dy\,\frac{dt}{t^{n}},
$$
where, as we will see, the limits exist in ${\mathcal S}'$. Consider the first integral
and, to simplify notation, define $\varphi(x)=x_k\omega(x)$.
Then, given $g\in{\mathcal S}$ we have to show that

$$
\int\left(\int_\eps^1\int \varphi\left(y+\frac{x-y}{t}\right)\, f(y)\, dy\,\frac{dt}{t^{n}}\right)\,g(x)\,dx
\to
\int\left(\int_0^1\int \varphi\left(y+\frac{x-y}{t}\right)\, f(y)\, dy\,\frac{dt}{t^{n}}\right)\,g(x)\,dx
$$
when $\eps\to 0$. It is enough to see that
\begin{equation}
\label{tiende a cero}
I_\eps:=\int_0^\eps\int\int \left|\varphi\left(y+\frac{x-y}{t}\right)\right|\,|f(y)|\,|g(x)|\,dx\, dy\,\frac{dt}{t^{n}}
\to 0.
\end{equation}
But, making the change of variable $z=\frac{x-y}{t}$ in the interior integral we have
$$
I_\eps=\int_0^\eps\int\int |\varphi(y+z)|\,|f(y)|\,|g(y+tz)|\,dz\, dy\,dt
\le \|g\|_{L^\infty(\R^n)}\|\varphi\|_{L^1(\R^n)}\|f\|_{L^1(\R^n)}\eps
$$
which proves (\ref{tiende a cero}). The integral defining $u_{k,2}$ can be treated
in the same way, indeed, defining now $\varphi(x)=\omega(x)$, the only difference
with the case of $u_{k,1}$ is the factor $y_k$ appearing in the integrand, but
it can be bounded assuming again that $f$ has compact support.

Now, for $\eps>0$ fixed, we can take the derivative inside the
integral, and therefore,

\begin{equation}
\label{derivada2}
\frac{\partial u_k}{\partial x_j}
=T_{kj,1}f+T_{kj,2}(y_kf)
\end{equation}
where $T_{kj,1}$ and $T_{kj,2}$ are of the form

\begin{equation}
\label{op gral}
Tf(x)= \lim_{\eps\to 0}\int_\eps^1\int
\frac{\partial}{\partial x_j}\left[\varphi\left(y+\frac{x-y}{t}\right)\right]
f(y)\, dy\,\frac{dt}{t^{n}}
\end{equation}
with $\varphi(x)=x_k\omega(x)$ for $T_{kj,1}$ and $\varphi(x)=\omega(x)$ for
$T_{kj,2}$.

We are going to prove continuity of operators of the form given
in (\ref{op gral}) where $\varphi\in C^\infty_0(B)$ with $B=B(0,\rho)$.
With this goal we decompose the operator as

$$
Tf=T_1f + T_2f
$$
where

\begin{equation}
\label{defT1}
T_1f(x)
=\lim_{\eps\to 0}\int_\eps^{\frac12}\int
\frac{\partial}{\partial x_j}\left[\varphi\left(y+\frac{x-y}{t}\right)\right]
f(y)\, dy\,\frac{dt}{t^{n}}
\end{equation}
and
$$
T_2f(x)=
\int_{\frac12}^1\int\frac{\partial}{\partial
x_j}\left[\varphi\left(y+\frac{x-y}{t}\right)\right]f(y)\, dy\,\frac{dt}{t^{n}}
$$

\medskip

An estimate of $\|T_2f\|_{L^2(\R^n)}$ for
$L^2$-functions $f$ vanishing outside $\O$ can be obtained easily
as we show in the following lemma.

\begin{lemma} If $f\in L^2(\R^n)$ vanishes outside $\O$ then
$$
\|T_2f\|_{L^2(\O)}
\le 2^{n}|\O|
\left\|\frac{\partial\varphi}{\partial x_j}\right\|_{L^\infty(\R^n)}
\|f\|_{L^2(\O)}
$$
\label{cota facil}
\end{lemma}

\proof
We have

\begin{equation}
\label{T2}
T_2f(x)= \int
\left\{\int_{\frac12}^1
\left(\frac{\partial\varphi}{\partial
x_j}\right)\left(y+\frac{x-y}{t}\right)\,\frac{dt}{t^{n+1}}\right\}
f(y)\, dy.
\end{equation}
Then,
$$
|T_2f(x)|
\le 2^n \left\|\frac{\partial\varphi}{\partial x_j}\right\|_{L^\infty(\R^n)}\int_\O |f(y)|\, dy
$$
and the result follows immediately using the Schwarz inequality. \qquad $\Box$

We now proceed to bound the operator $T_1$ in $L^2$. This will be done
using the Fourier transform. By standard density arguments it is enough
to bound the operator acting on $f$ smooth enough.
In the following lemma we give a simple form for $T_1$ in terms
of Fourier transforms.

\begin{lemma} For $f\in C_0^\infty(\R^n)$ we have
\begin{equation}
\label{op gral2}
\widehat{T_1f}(\xi)
=2\pi i\xi_j\int_0^{\frac12}  \widehat\varphi(t\xi)\widehat f((1-t)\xi) \,dt
\end{equation}
\label{fourier de T1}
\end{lemma}

\proof From (\ref{defT1}) we have

$$
T_1f=\lim_{\eps\to 0}T_{1,\eps}f,
$$
where
$$
T_{1,\eps}f(x)=\int_\eps^{\frac12}\int
\frac{\partial}{\partial x_j}\left[\varphi\left(y+\frac{x-y}{t}\right)\right]
f(y)\, dy\,\frac{dt}{t^{n}}
$$
and the limit is taken in ${\mathcal S}'$.

Now, we have

$$
\widehat{T_{1,\eps}f}(\xi)=
\int\int_\eps^{\frac12}\int
\frac{\partial}{\partial x_j}\left[\varphi\left(y+\frac{x-y}{t}\right)\right]
f(y)\,e^{-2\pi i x\cdot\xi} \,dy\,\frac{dt}{t^{n}}
 \, dx,
$$
and, since this triple integral exists, we can interchange the order of
integration. Therefore, integrating by parts we obtain

$$
\widehat{T_{1,\eps}f}(\xi)=2\pi i\xi_j
\int_\eps^{\frac12}\int \int
\varphi\left(y+\frac{x-y}{t}\right) f(y)
\,e^{-2\pi i x\cdot\xi} \, dx \,dy \,\frac{dt}{t^{n}},
$$
and making the change of variable
$$
z=y+\frac{(x-y)}{t}
$$
in the interior integral,

$$
\widehat{T_{1,\eps}f}(\xi)=2\pi i\xi_j
\int_\eps^{\frac12}\int \int
\varphi(z) \,e^{-2\pi i (tz+(1-t)y)\cdot\xi}
\,f(y)dz \,dy \,dt
$$
$$
=2\pi i\xi_j
\int_\eps^{\frac12}\int \int
\widehat\varphi(t\xi)
\,e^{-2\pi i(1-t)y\cdot\xi}
\,f(y) \,dy \,dt,
$$
and therefore,

$$
\widehat{T_{1,\eps}f}(\xi)
=2\pi i\xi_j\int_\eps^{\frac12}  \widehat\varphi(t\xi)\widehat f((1-t)\xi) \,dt,
$$
and taking $\eps\to 0$ we conclude the proof. \qquad $\Box$

Using the expression given in (\ref{op gral2}) we will give an estimate for the operator $T_1$ in $L^2$.
First we prove an auxiliary result.

\begin{lemma}
Define $C_{\varphi,\rho}=\rho^{-1}\|\varphi\|_{L^1(\R^n)}
+ \rho\left\|\frac{\partial^2\varphi}{\partial x_j^2}\right\|_{L^1(\R^n)}$.
Then,
$$
2\pi |\xi_j|\int_0^\infty |\widehat{\varphi}(t\xi)|\,dt\le
C_{\varphi,\rho}
$$
\label{cota1}
\end{lemma}
\proof
We have
$$
2\pi |\xi_j|\int_0^\infty |\widehat{\varphi}(t\xi)|\,dt
=2\pi |\xi_j|\int_0^{\frac1{2\pi\rho|\xi_j|}} |\widehat{\varphi}(t\xi)|\,dt
+2\pi |\xi_j|\int_{\frac1{2\pi\rho|\xi_j|}}^\infty |\widehat{\varphi}(t\xi)|\,dt
:= I + II
$$
Now,

$$
I\le \rho^{-1}\|\widehat\varphi\|_{L^\infty(\R^n)}\le
\rho^{-1}\|\varphi\|_{L^1(\R^n)}
$$
and
$$
II=2\pi\int_{\frac{1}{2\pi\rho|\xi_j|}}^\infty
\frac{ t^2|\xi_j|^2|\widehat\varphi(t\xi)|}{t^2|\xi_j|}\,dt
\le 2\pi\|\xi_j^2\widehat\varphi\|_{L^\infty(\R^n)}
\int_{\frac1{2\pi\rho|\xi_j|}}^\infty
\frac{1}{t^2|\xi_j|}\,dt
$$
but

$$
-4\pi^2\xi_j^2\widehat\varphi=\frac{\widehat{\partial^2\varphi}}{\partial x_j^2}
$$
and therefore,
$$
II\le\frac{1}{2\pi}
\left\|\frac{\widehat{\partial^2\varphi}}{\partial
x_j^2}\right\|_{L^\infty(\R^n)}
\int_{\frac1{2\pi\rho|\xi_j|}}^\infty
\frac{1}{t^2|\xi_j|}\,dt
= \rho\left\|\frac{\partial^2\varphi}{\partial x_j^2}\right\|_{L^1(\R^n)}
$$
and the lemma is proved. \qquad $\Box$

As a consequence of this lemma we obtain the following estimate for
the operator $T_1$.

\begin{lemma}
\label{cotaT1}
If $C_{\varphi,\rho}$ is the constant defined in the previous lemma, then
$$
\|T_1f\|_{L^2(\R^n)}
\le 2^{\frac{n-1}{2}} C_{\varphi,\rho}
\|f\|_{L^2(\R^n)}.
$$
\end{lemma}

\proof Applying the Schwarz inequality in (\ref{op gral2}) we have

$$
|\widehat{T_1f}(\xi)|^2
\le \left(\int_0^{\frac12} 2\pi |\xi_j| |\widehat\varphi(t\xi)|\,dt\right)
\left(\int_0^{\frac12} 2\pi |\xi_j| |\widehat\varphi(t\xi)||\widehat f((1-t)\xi)|^2
\,dt\right)
$$
and so, from Lemma \ref{cota1},

$$
|\widehat{T_1f}(\xi)|^2
\le C_{\varphi,\rho}\int_0^{\frac12} 2\pi |\xi_j|
|\widehat\varphi(t\xi)||\widehat
f((1-t)\xi)|^2 \,dt
$$
Then, integrating in $\xi$ and making the change of variable $\eta=(1-t)\xi$,
we obtain

$$
\int |\widehat{T_1f}(\xi)|^2\, d\xi
\le  C_{\varphi,\rho} \int_0^{\frac12}
\int \frac{2\pi}{(1-t)^{n+1}} |\eta_j|
\left|\widehat\varphi\left(\frac{t\eta}{1-t}\right)\right||\widehat f(\eta)|^2
\,
d\eta\,dt
$$
and, integrating first in the variable $t$ and making now the change
$s=t/(1-t)$, we get

$$
\int |\widehat{T_1f}(\xi)|^2\, d\xi
\le 2^{n-1} C_{\varphi,\rho}  \int\left(\int_0^1
2\pi |\eta_j| |\widehat\varphi(s\eta)|\,ds\right) \,|\widehat f(\eta)|^2 \,
d\eta
$$
therefore, applying again Lemma \ref{cota1},

$$
\int |\widehat{T_1f}(\xi)|^2\, d\xi
\le 2^{n-1} C^2_{\varphi,\rho}
\int
|\widehat f(\eta)|^2 \, d\eta
$$
and we conclude the proof recalling that the Fourier transform is an isometry in
$L^2(\R^n)$.\qquad $\Box$

\bigskip

Summing up we obtain the main result of this section.

\begin{theorem}
\label{main result}
If $T$ is the operator given in (\ref{op gral}) and $f$ vanishes outside $\O$,
then
$$
\|Tf\|_{L^2(\O)}\le C_{\varphi,\rho,\O}\|f\|_{L^2(\O)}
$$
\end{theorem}
with
$$
C_{\varphi,\rho,\O}=2^{\frac{n-1}{2}}\rho^{-1}\|\varphi\|_{L^1(\R^n)}
+ 2^{\frac{n-1}{2}}\rho\left\|\frac{\partial^2\varphi}{\partial x_j^2}\right\|_{L^1(\R^n)}
+ 2^{n}|\O|
\left\|\frac{\partial\varphi}{\partial x_j}\right\|_{L^\infty(\R^n)}.
$$
\proof The result follows immediately from Lemmas \ref{cota facil} and \ref{cotaT1}. \qquad $\Box$

\section{Dependence of the constant on $\O$}
\label{section2}
\setcounter{equation}{0}

An interesting question is what can be said, in terms of the geometry of the
domain $\O$, about the behavior of the
constant $C_{div,\O}$ in the estimate (\ref{estimate}).
Recall that we are assuming
that the domain $\O$ has diameter $R$ and that it is star-shaped with respect
to a ball of radius $\rho$ which, to simplify notation, we assume centered at
the origin.

It is known that the constant cannot be bounded independently of the ratio
$R/\rho$.
Indeed, this can be seen by the following elementary example which also
shows that, in some cases,

\begin{equation}
\label{cota inf}
C_{div,\O}\ge  c_1 (R/\rho)
\end{equation}
where $c_1$ is a constant independent of
$\O$.

Given positive numbers $a$ and $\ve$, consider the rectangular domain
$\O_{a,\ve}:= (-a,+a)\times (-\ve,\ve)$ and suppose that, for any
$f\in L^2(\O_{a,\ve})$ with vanishing mean value, there exists $\uu\in H_0^1(\O_{a,\ve})$ solving
(\ref{divu=f}) and satisfying the estimate (\ref{estimate}) with
a constant $C_{div,\O}=C_{a,\ve}$. Take $f(x_1,x_2)=x_1$ and the corresponding solution
$\uu$, then

$$
\|x_1\|^2_{L^2(\O_{a,\ve})}=\int_{\O_{a,\ve}} x_1 \d\uu
=-\int_{\O_{a,\ve}}u_1
=\int_{\O_{a,\ve}}x_2 \frac{\partial u_1}{\partial x_2}
$$
$$
\le \|x_2\|_{L^2(\O_{a,\ve})}
\left\|\frac{\partial u_1}{\partial x_2}\right\|_{L^2(\O_{a,\ve})}
\le C_{a,\ve} \|x_2\|_{L^2(\O_{a,\ve})}\|x_1\|_{L^2(\O_{a,\ve})}
$$
and so,

$$
\|x_1\|_{L^2(\O_{a,\ve})}
\le C_{a,\ve} \|x_2\|_{L^2(\O_{a,\ve})}
$$
but,

$$
\|x_1\|_{L^2(\O_{a,\ve})}=\frac2{\sqrt{3}}\ve^{\frac12} a^{\frac32}
\qquad \mbox{and} \qquad
\|x_2\|_{L^2(\O_{a,\ve})}=\frac2{\sqrt{3}}\ve^{\frac32} a^{\frac12}
$$
and therefore,

$$
C_{a,\ve}\ge (a/\ve)
$$
Consequently, if $a>\ve$, it follows that in this example
(\ref{cota inf}) holds.

\medskip

For the kind of domains that we are considering the following estimate for the constant
$C_{div,\O}$ is given in \cite{G}

$$
C_{div,\O} \le C_0 (R/\rho)^{n+1}
$$
with a constant $C_0$ independent of $\O$. The reader can check that
the result given in Theorem \ref{main result} recovers this estimate.
However, as we will show, this result can be improved.

Indeed, Theorem \ref{main result}
does not give a good estimate of the constant in terms of the function
$\varphi$ (or equivalently on $\rho$). Curiously, this is due to the estimate
obtained in Lemma \ref{cota facil} for the operator $T_2$ which in some sense is
easier to handle than $T_1$.
Then, in order to obtain a sharper bound,
we will give in the following lemmas a different argument to bound $T_2$.

\begin{lemma} If $1\le p<\frac{n}{n-1}$ then,
$$
\|T_2f\|_{L^p(\R^n)}
\le \frac{2^{\frac{n}{p'}}}{(1-\frac{n}{p'})}\,
\left\|\frac{\partial\varphi}{\partial x_j}\right\|_{L^1(\R^n)}
\|f\|_{L^p(\R^n)}
$$
\label{cota dificil}
\end{lemma}

\proof From (\ref{T2}) we have

$$
|T_2f(x)|
\le \int_{\frac12}^1 \int
\left|\left(\frac{\partial\varphi}{\partial
x_j}\right)\left(y+\frac{x-y}{t}\right)\right|
|f(y)|\, dy \,\frac{dt}{t^{n+1}}
$$
Making the change of variable
$$
z=y+\frac{x-y}{t}
$$
in the interior integral, we obtain

$$
|T_2f(x)|
\le 2 \int_{\frac12}^1 \int
\left|\frac{\partial\varphi}{\partial x_j}(z)\right|
\left|f\left(\frac{tz-x}{t-1}\right)\right| \frac{1}{(1-t)^n}\, dz \,dt .
$$
Applying now the Minkowski inequality for integrals we have

$$
\|T_2f\|_{L^p(\R^n)}
\le 2 \int_{\frac12}^1 \int
\left|\frac{\partial\varphi}{\partial x_j}(z)\right|
\left(\int \left|f\left(\frac{tz-x}{t-1}\right)\right|^pdx\right)^\frac{1}{p}
\frac{1}{(1-t)^n}\, dz \,dt
$$
and, by the change of variable
$$
\overline x=\frac{tz-x}{t-1}
$$
in the interior integral, it
follows that

$$
\|T_2f\|_{L^p(\R^n)}
\le 2 \left\|\frac{\partial\varphi}{\partial x_j}\right\|_{L^1(\R^n)}
\|f\|_{L^p(\R^n)}
\int_{\frac12}^1  \frac{1}{(1-t)^{\frac{n}{p'}}}\,dt
$$
therefore, since $p'>n$, the integral on the right hand side of this
inequality is finite and so we obtain the lemma. \qquad $\Box$

Unfortunately the restriction for the value of $p$ in the previous
lemma excludes the case $p=2$. However, using well known interpolation
theorems we can obtain an estimate for the $L^2$ case.

\begin{lemma}
\label{interpolation}
If $f\in L^2(\R^n)$ vanishes outside $\O$ then, for $1\le p<\frac{n}{n-1}$,

$$
\|T_2f\|_{L^2(\O)}
\le  \frac{2^{\frac{n}{2}}}{(1-\frac{n}{p'})^{\frac{p}{2}}} |\O|^{1-\frac{p}{2}}\left\|\frac{\partial\varphi}{\partial
x_j}\right\|^{\frac{p}{2}}_{L^1(\R^n)}
\left\|\frac{\partial\varphi}{\partial x_j}\right\|^{1-\frac{p}{2}}_{L^\infty(\R^n)}
\|f\|_{L^2(\O)}
$$
\end{lemma}

\proof From the definition of $T_2$ (\ref{T2}) it is easy to see that
$$
\|T_2f\|_{L^\infty(\O)}
\le 2^{n}|\O|
\left\|\frac{\partial\varphi}{\partial x_j}\right\|_{L^\infty(\R^n)}
\|f\|_{L^\infty(\O)}.
$$
Then, the result follows immediately from this estimate together
with Lemma \ref{cota dificil} and the
well known interpolation inequality
$$
\|T_2\|_{{\mathcal L}(L^2,L^2)}
\le \|T_2\|^{\frac{p}{2}}_{{\mathcal L}(L^p,L^p)}
\|T_2\|^{1-\frac{p}{2}}_{{\mathcal L}(L^\infty,L^\infty)}.
\qquad \Box
$$

Summing up we obtain the following estimate in terms of the function $\varphi$.

\begin{theorem}
\label{main result2}
If $T$ is the operator given in (\ref{op gral}), $f$ vanishes outside $\O$,
and $1\le p<\frac{n}{n-1}$, then

$$
\|Tf\|_{L^2(\O)}\le \left(2^{\frac{n-1}{2}} C_{\varphi,\rho}
+ \frac{2^{\frac{n}{2}}}{(1-\frac{n}{p'})^{\frac{p}{2}}}
\widetilde C_{\varphi,p}|\O|^{1-\frac{p}{2}}\right) \|f\|_{L^2(\O)}
$$
where
$$
C_{\varphi,\rho}=\rho^{-1}\|\varphi\|_{L^1(\R^n)}
+ \rho\left\|\frac{\partial^2\varphi}{\partial x_j^2}\right\|_{L^1(\R^n)}
$$
and
$$
\widetilde C_{\varphi,p}
=\left\|\frac{\partial\varphi}{\partial x_j}\right\|^{\frac{p}{2}}_{L^1(\R^n)}
\left\|\frac{\partial\varphi}{\partial
x_j}\right\|^{1-\frac{p}{2}}_{L^\infty(\O)}
$$

\end{theorem}
\proof The result follows immediately from Lemmas \ref{cotaT1} and \ref{interpolation}. \qquad $\Box$

We want to bound $\|\frac{\partial u_k}{\partial x_j}\|_{L^2(\R^n)}$ using the expression (\ref{derivada2}).
This is the goal of the following theorem.

In what follows $C_n$ denotes a constant depending only on $n$,
not necessarily the same at each occurrence, and
$A\sim B$ means that $A/B$ is bounded by above and below by positive constants
which may depend on $n$ and $p$ only.

\begin{theorem}
\label{main result3}
Let $\O\subset\R^n$ be a bounded domain of diameter $R$ which is star-shaped
with respect to
a ball $B\subset\O$ of radius $\rho$ and $\uu\in H^1_0(\O)$ be
the solution of (\ref{divu=f}) given by (\ref{solucion}). Then,
there exists a constant $C_n$ such that

$$
\left\|D\uu\right\|_{L^2(\O)}
\le C_n
\,\frac{R}{\rho}\, \left(\frac{|\O|}{|B|}\right)^{\frac{n-2}{2(n-1)}}
\left(\log \frac{|\O|}{|B|}\right)^{\frac{n}{2(n-1)}}
\|f\|_{L^2(\O)}
$$

\end{theorem}

\proof
As we have mentioned, both operators
on the right hand side of (\ref{derivada2}) are of the form given in (\ref{op gral}).
We will estimate the term $T_{kj,2}(y_kf)$ which is the worst part due to the
presence of $y_k$. The reader can check that the term $T_{kj,1}f$ can be
bounded analogously.

For $T_{kj,2}$ the function $\varphi$ is exactly
$\omega$, which is supported in $B(0,\rho)$ and has integral equal to
one. Therefore, $\varphi$ can be taken as

$$
\varphi(x)=\rho^{-n}\psi(\rho^{-1}x),
$$
where $\psi$ is a smooth function supported in the unit ball and
with integral equal to one. Then,

$$
\frac{\partial\varphi}{\partial x_j}(x)
=\rho^{-n-1}\frac{\partial\psi}{\partial x_j}(\rho^{-1}x)
\qquad , \qquad
\frac{\partial^2\varphi}{\partial x^2_j}(x)
=\rho^{-n-2}\frac{\partial^2\psi}{\partial x^2_j}(\rho^{-1}x)
$$
and so,
\begin{equation}
\label{cte1}
C_{\varphi,\rho}=\rho^{-1}\|\varphi\|_{L^1(\R^n)}
+ \rho\left\|\frac{\partial^2\varphi}{\partial x_j^2}\right\|_{L^1(\R^n)}
\sim \, \rho^{-1}
\end{equation}
and
\begin{equation}
\label{cte2}
\widetilde C_{\varphi,p}
=\left\|\frac{\partial\varphi}{\partial x_j}\right\|^{\frac{p}{2}}_{L^1(\R^n)}
\left\|\frac{\partial\varphi}{\partial
x_j}\right\|^{1-\frac{p}{2}}_{L^\infty(\R^n)}
\sim \, \rho^{-1-n(1-\frac{p}{2})}.
\end{equation}

Therefore, applying Theorem \ref{main result2} for $T=T_{kj,2}$, using $|y_k|\le R$
and the relations (\ref{cte1})
and (\ref{cte2}), we obtain, for $1\le p<\frac{n}{n-1}$,

\begin{align*}
\left\|D\uu\right\|_{L^2(\O)}
&\le C_n
\,\frac{R}{\rho}\,
\frac{1}{(1-\frac{n}{p'})^{\frac{p}{2}}}
\left(\frac{|\O|}{|B|}\right)^{1-\frac{p}{2}}
\|f\|_{L^2(\O)}\\
{}& = C_n
\,\frac{R}{\rho}\,
\frac{1}{(1-\frac{n}{p'})^{\frac{p}{2}}}
\left(\frac{|\O|}{|B|}\right)^{\frac{n-2}{2(n-1)}}
\left(\frac{|\O|}{|B|}\right)^{\frac12\left(\frac{n}{n-1}-p\right)}
\|f\|_{L^2(\O)}
\end{align*}

Now, assuming that $\frac{|\O|}{|B|}$ is large enough, we can choose $p$
such that
$$
\frac12\left(\frac{n}{n-1}-p\right)
=\frac{1}{\log \frac{|\O|}{|B|}}
$$
obtaining
$$
\left\|D\uu\right\|_{L^2(\O)}
\le C_n
\,\frac{R}{\rho}\,
\frac{1}{(1-\frac{n}{p'})^{\frac{p}{2}}}
\left(\frac{|\O|}{|B|}\right)^{\frac{n-2}{2(n-1)}}
e\|f\|_{L^2(\O)},
$$
and so, we conclude the proof using that
$$
1-\frac{n}{p'}=\frac{(n-1)}{p}\left(\frac{n}{n-1}-p\right)
=\frac{2(n-1)}{p\log\left(\frac{|\O|}{|B|}\right)}
$$
and $p<\frac{n}{n-1}$. \qquad \qquad $\Box$

\begin{remark} In the particular case $n=2$ the theorem gives

$$
\left\|D\uu\right\|_{L^2(\O)}
\le C
\,\left(\frac{R}{\rho}\right)
\,\log \left(\frac{R}{\rho}\right)
\|f\|_{L^2(\O)}
$$

In view of the example given above this estimate is almost optimal
(i.e., optimal up to the logarithmic factor).
\end{remark}

\section{The constant in the Korn inequality}
\label{section3}
\setcounter{equation}{0}

As it is well known, Korn type inequalities are strongly connected with the
existence of solutions of (\ref{divu=f}) satisfying (\ref{estimate}).
For example, in the particular case of two dimensional simple connected domains with a $C^1$ boundary,
the explicit relation between the best constant in (\ref{estimate}) and that
in the so-called second case of Korn inequality was given in \cite{HP}.
More generally, for arbitrary domains in $n$ dimensions, $n\ge 2$,
the Korn inequality can be derived from the existence of solutions of
the divergence satisfying (\ref{estimate}), and therefore, information on
the constant in the Korn inequality can be obtained from estimates for the
constant in (\ref{estimate}).

A lot of work has been done in order to obtain the behavior of the
constant in the different versions of Korn inequality in terms of
the domain (see \cite{H} and its references).

We are going to show how our results in the previous section
can be used to obtain estimates for the constant in the second
case of Korn inequality. Let us mention that domains which are star-shaped
with respect to a ball were considered by Kondratiev and Oleinik in \cite{KO1,KO2}
where the authors obtain sharp estimates for the constant in a Korn inequality in
terms of $R/\rho$.
However, their results are for a different type of Korn inequality than the
one that we are considering and it is not clear what is the relation
between the constants in the two different Korn type inequalities.

For a vector field $\vv\in H^1(\O)^n$,
$\ve(\vv)$ and $\mu(\vv)$ denote its symmetric and skew symmetric part
respectively, i. e.,
$$
\ve_{ij}(\vv)=\frac12 \left(\frac{\partial v_i}{\partial x_j}
+ \frac{\partial v_j}{\partial x_i}\right)
$$
and
$$
\mu_{ij}(\vv)
=\frac12 \left(\frac{\partial v_i}{\partial x_j}
- \frac{\partial v_j}{\partial x_i}\right)
$$

Then, the so-called second case of Korn inequality states that there
exists a constant $C_{K,\O}$ such that
$$
\|D\vv\|_{L^2(\O)}
\le C_{K,\O}\|\ve(\vv)\|_{L^2(\O)}
$$
for vector fields $\vv\in H^1(\O)^n$ satisfying
\begin{equation}
\label{caso2}
\int_\O \mu_{ij}(\vv)= 0\, ,
\qquad \mbox{for\,}\quad i,j=1,\ldots, n.
\end{equation}

The argument used in the proof of the following theorem is known
but we include it for the sake of completeness. For an arbitrary
domain $\O$ we will say that it admits a right inverse of the divergence
with constant $C_{div,\O}$ if for any $f\in L^2(\O)$, such that $\int_\O f=0$,
there exists $\uu\in H^1_0(\O)^n$ satisfying

$$
\d\uu=f
$$
and
$$
\left\|D\uu\right\|_{L^2(\O)}
\le C_{div,\O} \|f\|_{L^2(\O)}.
$$

\begin{theorem}
\label{div implica Korn}
If $\O$ admits a right inverse of the divergence
with constant $C_{div,\O}$, then the second case of Korn inequality
holds in $\O$ with a constant $C_{K,\O}$ which satisfies
$$
C_{K,\O}\le (1+4n^2)^{1/2} C_{div,\O}
$$
\end{theorem}

\proof
Let $\vv\in H^1(\O)^n$ such that (\ref{caso2}) holds. By density
we can assume that $\vv$ is smooth.
By orthogonality we have

$$
\|D\vv\|^2_{L^2(\O)}
=\|\ve(\vv)\|^2_{L^2(\O)}+\|\mu(\vv)\|^2_{L^2(\O)}
$$
and so, observing that $C_{div,\O}\ge 1$, it is enough to prove
that
\begin{equation}
\label{comparacion}
\|\mu(\vv)\|^2_{L^2(\O)}\le 4 n^2 C^2_{div,\O}\|\ve(\vv)\|^2_{L^2(\O)}.
\end{equation}
Given $i$ and $j$, since
$$
\int_\O \mu_{ij}(\vv) = 0,
$$
there exists $\uu^{ij}\in H_0^1(\O)^n$ such that
$$
\d\uu^{ij}=\mu_{ij}(\vv)
$$
and
\begin{equation}
\label{ri}
\left\|D\uu^{ij}\right\|_{L^2(\O)}
\le C_{div,\O} \|\mu_{ij}(\vv)\|_{L^2(\O)}.
\end{equation}
Then,
$$
\|\mu_{ij}(\vv)\|^2_{L^2(\O)}
=\int_\O \mu_{ij}(\vv)\,\d\uu^{ij}
=-\int_\O \nabla\mu_{ij}(\vv)\cdot\uu^{ij}
$$
but,
$$
\frac{\partial\mu_{ij}(\vv)}{\partial x_k}
=\left(\frac{\partial \ve_{ik}(\vv)}{\partial x_j}
- \frac{\partial \ve_{jk}(\vv)}{\partial x_i}\right)
$$
and so,
$$
\|\mu_{ij}(\vv)\|^2_{L^2(\O)}
=-\sum_{k=1}^n
\int_\O\left(\frac{\partial \ve_{ik}(\vv)}{\partial x_j}
- \frac{\partial \ve_{jk}(\vv)}{\partial x_i}\right)u^{ij}_k
=\sum_{k=1}^n\int_\O
\left(\ve_{ik}(\vv)\frac{\partial u^{ij}_k}{\partial x_j}
-\ve_{jk}(\vv)\frac{\partial u^{ij}_k}{\partial x_i}\right),
$$
and using now (\ref{ri}) we obtain,
$$
\|\mu_{ij}(\vv)\|^2_{L^2(\O)}
\le C_{div,\O} \|\mu_{ij}(\vv)\|_{L^2(\O)}
\sum_{k=1}^n \left(\|\ve_{ik}(\vv)\|_{L^2(\O)} + \|\ve_{jk}(\vv)\|_{L^2(\O)}\right).
$$
Therefore,
$$
\|\mu_{ij}(\vv)\|_{L^2(\O)}
\le C_{div,\O} \, n^{1/2}
\left\{\sum_{k=1}^n \left(\|\ve_{ik}(\vv)\|_{L^2(\O)} + \|\ve_{jk}(\vv)\|_{L^2(\O)}\right)^2\right\}^{1/2}
$$
and then
$$
\|\mu_{ij}(\vv)\|^2_{L^2(\O)}
\le 2 C_{div,\O}^2 \, n
\sum_{k=1}^n \left(\|\ve_{ik}(\vv)\|^2_{L^2(\O)} + \|\ve_{jk}(\vv)\|^2_{L^2(\O)}\right).
$$
Finally, summing now in $i$ and $j$ we obtain
(\ref{comparacion}). \qquad $\Box$

Consequently, using the results of the previous section we obtain an estimate for the Korn inequality
in star-shaped domains.

\begin{theorem}
\label{main result Korn}
Let $\O\subset\R^n$ be a bounded domain of diameter $R$ which is star-shaped
with respect to a ball $B\subset\O$ of radius $\rho$. Then,
there exists a constant $C_n$ such that, for
all $\vv\in H^1(\O)^n$ satisfying
$\int_\O \mu_{ij}(\vv)=0$, for $i,j=1,\ldots, n$,

$$
\|D\vv\|_{L^2(\O)}
\le C_n
\,\frac{R}{\rho}\, \left(\frac{|\O|}{|B|}\right)^{\frac{n-2}{2(n-1)}}
\left(\log \frac{|\O|}{|B|}\right)^{\frac{n}{2(n-1)}}
\|\ve(\vv)\|_{L^2(\O)}
$$
\end{theorem}

\proof The result follows immediately from Theorems \ref{main result2} and
\ref{div implica Korn}. \qquad $\Box$

\section{The constant in the improved Poincar\'e inequality}
\label{section4}
\setcounter{equation}{0}

In this section we consider another well known inequality which is related
with the existence of right inverses of the divergence, namely, the so-called
improved Poincar\'e inequality. To recall this inequality we need to introduce
some notation. For a bounded domain $\O\subset\R^n$ and any $x\in\O$ we denote with
$d(x)$ the distance of $x$ to the boundary of $\O$.
Then, the improved
Poincar\'e inequality states that there exists a constant $C_{iP,\O}$
such that, for any $f\in H^1(\O)\cap L^2_0(\O)$, where
$L^2_0(\O):=\{f\in L^2(\O)\,:\,\int_\O f=0\}$,

\begin{equation}
\label{iP}
\|f\|_{L^2(\O)}\le C_{iP,\O}\|d\nabla f\|_{L^2(\O)}.
\end{equation}

It is known that this inequality is valid for Lipschitz domains and,
more generally, for John domains (see for example \cite{BS,DD,HS}).

For the star-shaped domains that we are considering in this paper, the argument
given in \cite{DD}, applied in this particular case, can be used
to show that
\begin{equation}
\label{cota fusco}
C_{iP,\O}\le C_n (R/\rho)^{n+1},
\end{equation}
indeed, this was done in \cite[Prop. 5.2]{BCF} for the analogous inequality
in $L^1$, but it is easy to see that the arguments extend straightforward to the
$L^2$ case. We are going to show that the dependence on $R/\rho$ can be improved
using our estimates of Section \ref{section2}, at least in the two dimensional case.

Recently, in \cite{DMRT} the relation between Poincar\'e type inequalities
and solutions of the divergence was analyzed in a very general context.
A particular case of the results in that paper says that
the improved Poincar\'e inequality (\ref{iP})
implies the existence of a right inverse of the divergence as an operator from
$H^1_0(\O)^n$ to $L^2_0(\O)$.
Let us reproduce the argument given in that paper in this particular case
for the sake of completeness. With this purpose we need to use a
Whitney decomposition of $\O$, i. e., a sequence of cubes $\{Q_j\}$
with pairwise disjoints interiors and such that, if $d_j$ and $\ell_j$
are the distance of $Q_j$ to the boundary of $\O$ and the length of
its edges respectively, then $d_j/\ell_j$ is bounded by above and below
by positive constants depending only on $n$.
Associated with this decomposition there is a
partition of unity $\{\phi_j\}$, namely, $\sum_j \phi_j=1$ with $\phi_j\in C_0^\infty(\tQ_j)$
where $\tQ_j$ is an expansion of $Q_j$ still with diameter
proportional to its distance to the boundary of $\O$ (see for example \cite{S} for details).

\begin{lemma}
\label{descomposicion}
If the improved Poincar\'e inequality (\ref{iP}) is satisfied in
$\O$ then, given $f\in L_0^2(\O)$ and a Whitney decomposition of $\O$, there exists
a sequence $\{f_j\}$ such that $f_j\in L_0^2(\tQ_j)$, $f=\sum_j f_j$, and
\begin{equation}
\label{desc1}
\|f\|^2_{L^2(\O)}\le C_n \sum_j\|f_j\|^2_{L^2(\tQ_j)}
\end{equation}
and
\begin{equation}
\label{desc2}
\sum_j\|f_j\|^2_{L^2(\tQ_j)}\le C_n (1+C_{iP,\O}) \|f\|^2_{L^2(\O)}.
\end{equation}
\end{lemma}

\proof
First we observe that, by duality, (\ref{iP}) implies
that, for all $f\in L_0^2(\O)$, there exists $\vv\in L^2(\O)^n$ such that
\begin{equation}
\label{dual}
\d\vv=f \quad \mbox{in}\ \O \ , \ \vv\cdot\nn=0 \quad \mbox{on}\ \partial\O
\end{equation}
and
\begin{equation}
\label{dual2}
\left\|\frac{\vv}{d}\right\|_{L^2(\O)}\le C_{iP,\O} \|f\|_{L^2(\O)},
\end{equation}
where both equations in (\ref{dual}) has to be understood in a distributional
sense.

Indeed,
$$
L(\nabla g)=\int_\O fg
$$
defines a linear form on the subspace of $L^2(\O)^n$ formed
by the gradient vector fields. $L$ is well defined because
$\int_\O f=0$. Moreover, it follows from (\ref{iP}) that
$$
|L(\nabla g)|=\left|\int_\O f (g-\overline g)\right|
\le C_{iP,\O} \|f\|_{L^2(\O)}\|d\nabla g\|_{L^2(\O)}
$$
where $\overline g$ is the average of $g$ in $\O$.

By the Hahn-Banach theorem $L$ can be extended as a linear continuous functional
to the space $L_d^2(\O)^n$, where $L^2_d(\O)$ denotes the Hilbert space with norm
$\|f\|_{L_d^2}:=\|df\|_{L^2}$,
and therefore, there exists $\vv\in L^2_{d^{-1}}(\O)^n$
satisfying (\ref{dual2}) and such that
$$
L(\ww)=\int_\O\vv\cdot\ww \quad  \quad
\forall\,\ww\in L_d^2(\O)^n,
$$
in particular,
$$
\int_\O\vv\cdot\nabla g=\int fg \quad \forall g\in H^1(\O)
$$
which is equivalent to (\ref{dual}).

Given now $f\in L_0^2(\O)$ let $\vv\in L^2(\O)^n$ satisfying (\ref{dual})
and (\ref{dual2}) and define
$$
f_j=\d(\phi_j\vv).
$$
Then, we have
$$
f=\d\vv=\d\Big(\vv \sum_{j}\phi_j\Big)=\sum_{j}\d(\phi_j\vv)=\sum_j f_j.
$$
Since ${\rm supp\,}\phi_j\subset\tQ_j$ we have
${\rm supp\,}f_j\subset\tQ_j$ and $\int f_j=0$.

Moreover, using the finite superposition (with constant depending only on $n$)
of the expanded cubes $\tQ_j$, we obtain immediately (\ref{desc1}).
On the other hand, using again the finite superposition
and that $\|\phi_j\|_{L^\infty}\le 1$ and $\|\nabla\phi_j\|_{L^\infty}\le C/d_j$,
we have
$$
\|f_j\|^2_{L^2(\tQ_j)}
\le C_n \left\{ \|f\|^2_{L^2(\tQ_j)} + \left\|\frac{\vv}{d}\right\|^2_{L^2(\tQ_j)}\right\},
$$
and therefore, (\ref{desc2}) follows from (\ref{dual2}).\qquad $\Box$

\begin{theorem}
\label{improved implica div}
If the improved Poincar\'e inequality (\ref{iP}) is satisfied in
$\O$ then, $\O$ admits a right inverse of the divergence
with constant $C_{div,\O}$ which satisfies
\begin{equation}
\label{comparacion2}
C_{div,\O}\le C_n (1+C_{iP,\O}).
\end{equation}
\end{theorem}

\proof Given $f\in L_0^2(\O)$ let $f_j$ be the functions given in
the previous lemma. Since $f_j\in L_0^2(\tQ_j)$, there exists $\uu_j\in H_0^1(\tQ_j)^n$
such that
$$
\d\uu_j=f_j \quad \mbox{and} \quad
\|D\uu_j\|_{L^2(\tQ_j)}\le C_n \|f_j\|_{L^2(\tQ_j)},
$$
indeed, a scaling argument shows that the constant in
this inequality is independent of the size of the cube.
Then, $\uu=\sum_j\uu_j\in H_0^1(\O)^n$ is a solution of
$\d\uu=f$. Moreover, it follows from (\ref{desc2}) that
$$
\|D\uu\|_{L^2(\O)}\le C_n (1+C_{iP,\O})\|f\|_{L^2(\O)}
$$
and the theorem is proved.\qquad $\Box$

In view of the previous theorem a natural question is whether
the converse is also true. To the author knowledge this is not
known. However, a weaker result will allow us to obtain an estimate
for the constant in the improved Poincar\'e inequality for planar star-shaped domains.
In fact, we will see that the converse can be proved if we assume
that the following inequality is satisfied in $\O$ (actually this is one
of the many results called ``Hardy inequality'' although, at least to the author
knowledge, Hardy proved only the one dimensional case).

\begin{equation}
\label{hardy}
\left\|\frac{g}{d}\right\|_{L^2(\O)}\le C_{H,\O}\|\nabla g \|_{L^2(\O)}
\qquad \forall\, g\in H^1_0(\O).
\end{equation}

It is known that this inequality is valid for a very large class of domains
(see for example \cite{H2,KK,M,OK}).

\begin{theorem}
\label{div + hardy implica improved}
If $\O$ admits a right inverse of the divergence
with constant $C_{div,\O}$ and the Hardy inequality (\ref{hardy})
is satisfied in $\O$ then,
the improved Poincar\'e inequality (\ref{iP}) is valid in
$\O$ with a constant $C_{iP,\O}$ such that
\begin{equation}
\label{comparacion3}
C_{iP,\O}
\le C_{H,\O} C_{div,\O}.
\end{equation}
\end{theorem}

\proof Given $f\in L_0^2(\O)$ let $\uu\in H_0^1(\O)^n$ be
such that
\begin{equation}
\label{div2}
\d\uu=f \quad \mbox{and} \quad
\|D\uu\|_{L^2(\O)}\le C_{div,\O} \|f\|_{L^2(\O)}.
\end{equation}
Then,
$$
\|f\|^2_{L^2(\O)}
=\int_\O f \d\uu=-\int_\O \nabla f\cdot\uu
\le \|d\nabla f\|_{L^2(\O)} \left\|\frac{\uu}{d}\right\|_{L^2(\O)}
\le C_{H,\O}\|d\nabla f\|_{L^2(\O)} \|D\uu\|_{L^2(\O)}
$$
and using (\ref{div2}) we conclude the proof.\qquad $\Box$

In order to apply this theorem together with our results of
Section \ref{section2} we need to know estimates for $C_{H,\O}$.
For example, for simply connected (in particular for star-shaped)
planar domains it has been proved that
\begin{equation}
\label{constante de hardy}
C_{H,\O}\le 4,
\end{equation}
see \cite{A,Ba}.

Therefore, using this estimate and the results of Section \ref{section2},
we obtain an estimate for the constant $C_{iP}$ which improves (\ref{cota fusco}).

\begin{theorem}
\label{main result IP}
Let $\O\subset\R^2$ be a bounded domain of diameter $R$ which is star-shaped
with respect to a ball $B\subset\O$ of radius $\rho$. Then, there exists a positive constant $C$
such that, for all $f\in H^1(\O)\cap L^2_0(\O)$, we have
$$
\|f\|_{L^2(\O)}
\le C \,\left(\frac{R}{\rho}\right)
\,\log \left(\frac{R}{\rho}\right)\|d\nabla f\|_{L^2(\O)}.
$$
\end{theorem}

\proof The result follows immediately from Theorems \ref{main result3} and
\ref{div + hardy implica improved} and inequality (\ref{constante de hardy}). \qquad $\Box$

To finish the paper let us mention that the bound given in the previous theorem
is almost optimal. Indeed, in view of Theorem \ref{improved implica div}, the
same example given in Section \ref{section2} shows that in some cases
$C_{iP}\ge c_1(R/\rho)$, where $c_1$ is a constant.


\begin{thebibliography}{99}



\bibitem{ADM} G. Acosta, R. G. Dur\'an, and M.A. Muschietti,
{\sl Solutions of the divergence operator on John domains},
Adv. Math.  206(2) , 373-401, 2006.

\bibitem{A} A. Ancona,
{\sl On strong barriers and an inequality of hardy
for domains in $\R^n$},
J. London
Math. Soc. (2)34, 274-290, 1986.

\bibitem{Ba} R. Ba\~nuelos,
{\sl Four unknown constants},
2009 Oberwolfach workshop on Low Eigenvalues of Laplace and Schr\"odinger Operators,
Oberwolfach reports 6(1), 2009.


\bibitem{BCF} M. Barchiesi, F. Cagnetti, and N. Fusco,
{\sl Stability of the steiner symmetrization of convex sets},
preprint, CNA, Carnegie Mellon University, USA, 2009.


\bibitem{BS} H. B. Boas, and E. J. Straube,
{\sl Integral inequalities of Hardy and Poincare type},
Proc. Amer. Math. Soc. 103(1), 172-176, 1988.

\bibitem{B} M. E. Bogovskii,
{\sl Solution of the first boundary value problem for the
equation of continuity of an incompressible medium},
Soviet Math. Dokl. 20, 1094-1098, 1979.


\bibitem{CZ} A. P. Calder\'on, and A. Zygmund,
{\sl On singular integrals},
Amer. J. Math. 78(2), 289-309, 1956.

\bibitem{CM} M. Costabel, and A. McIntosh,
{\sl On Bogovskii and regularized Poincar\'e
integral operators for de Rham complexes on Lipschitz domains},
Math. Z. 265(2) , 297–320, 2010.

\bibitem{DD} I. Drelichman, and R. G. Dur\'an,
{\sl Improved Poincar\'e inequalities with weights},
J. Math. Anal. Appl. 347(1), 286–293, 2008.

\bibitem{DM} R. G. Dur\'an, and M.A. Muschietti,
{\sl An explicit right inverse of the divergence
operator which is continuous in weighted norms},
Studia Mathematica  148(3), 207-219, 2001.

\bibitem{DMRT} R. G. Dur\'an, M.A. Muschietti, E. Russ, and P. Tchamitchian,
{\sl Divergence operator and Poincar\'e inequalities on arbitrary bounded domains},
Complex Var. Elliptic Equ. 55(8-10) ), 795–816, 2010.

\bibitem{G} G.P. Galdi, An Introduction to the Mathematical Theory of the Navier-
Stokes Equations. Vol. I. Linearized Steady Problems, Springer Tracts in
Natural Philosophy 38, Springer-Verlag 1994.

\bibitem{GHH} M. Geissert, H. Heck, and M. Hieber,
{\sl On the equation $div u=g$ and Bogovskii's operator
in Sobolev spaces of negative order},
Operator Theory: Advances and Applications 168, 113-121, 2006.

\bibitem{H2} P. Hajlasz,
{\sl Pointwise Hardy inequalities},
Proc. Amer. Math. Soc. 127(2), 417–423, 1999.

\bibitem{H} C. O. Horgan,
{\sl Korn's inequalities and their applications in continuum mechanics},
SIAM Review 37(4), 491-511, 1995.

\bibitem{HP} C. O. Horgan, and L. E. Payne,
{\sl On inequalities of Korn, Friedrichs and Babuška-Aziz}, Arch.
Rational Mech. Anal. 82(2), 165-179, 1983.

\bibitem{HS} R. Hurri-Syrj\"anen,
{\sl An improved Poincar\'e inequality},
Proc. Amer. Math. Soc. 120(1), 213-222, 1994.

\bibitem{KK} J. Kinnunen, and R. Korte,
{\sl Characterizations for the Hardy inequality},
Around the research of Vladimir Maz'ya. I, 239–254,
Int. Math. Ser. (N. Y.), 11, Springer, New York, 2010.

\bibitem{KO1} V. A. Kondratiev and O. A. Oleinik,
{\sl On Korn's inequalities}, C. R. Acad. Sci. Paris 308(16), 483-487, 1989.

\bibitem{KO2} V. A. Kondratiev, and O. A. Oleinik,
{\sl Hardy's and Korn's type inequalities and their applications},
Rend. Mat. Appl. Serie VII 10(3), 641-666, 1990.

\bibitem{M} V. G. Maz'ya, Sobolev spaces, Springer-Verlag, 1985.

\bibitem{OK} B. Opic, and A. Kufner,
Hardy-type inequalities, Research Notes in Mathematics 219, Pitman,
1990.

\bibitem{S} E. M. Stein, Singular integrals and
differentiability properties of functions, Princeton Univ.
Press, 1970.

\end{thebibliography}
\end{document}